\documentclass{article}
\usepackage{longtable}
\usepackage{multirow}
\usepackage{scalefnt}
\usepackage{amsmath,amsthm,amsfonts,amssymb}
\usepackage{bm}
\usepackage{here}
\usepackage{caption}
\usepackage[margin=10truemm]{geometry}
\numberwithin{equation}{section}
\title{\huge On the diophantine equation $X^6 - Y^6 = W^n - Z^n,\ n=2,3,4$}
\author{Seiji Tomita }
\date{}
\begin{document}
\maketitle
\begin{abstract}
In this paper, we proved that there are infinitely many integer solutions of $X^6 - Y^6 = W^n - Z^n,\ n=2,3,4$.

\end{abstract}
\vskip3\baselineskip

\section{ Introduction}

The diophantine equation $x^n+y^n=z^n+w^n, n=2,3,4$ have been considered by many mathematicians.
Euler\cite{a} showed the parametric solution for $n=4$.
Choudhry\cite{b} considered the equation  $x^4+hy^4=z^4+hw^4$.
Izadi\cite{c} and Nabardi showed that there are infinitely many integer solutions of $x^4+y^4=2(z^4+w^4)$.
Janfada\cite{d} and Nabardi considered the equation $x^4+y^4=h(z^4+w^4)$.
Shabani-Solt\cite{f},  Yusefnejad, and  Janfada considered the equation $x^6+ky^3=z^6+kw^3$ with $k<500$.
Muthuvel\cite{e} and Venkatraman showed that there are infinitely many integer solutions of $X^4-Y^4=R^2-S^2$.

The objective of this work is to prove that there are infinitely many integer solutions of $X^6 - Y^6 = W^n - Z^n,\ n=2,3,4.$

\vskip3\baselineskip

\section{Solving the diophantine equation $X^6 - Y^6 = W^2 - Z^2$}
\begin{equation}
X^6 - Y^6 = W^2 - Z^2
\end{equation}

We give three methods to prove that there are infinitely many integer solutions of $X^6 - Y^6 = W^2 - Z^2$.
\vskip3\baselineskip
2.1. Method-1
\vskip\baselineskip

 Let $X=a$ and $Y=b$ where $a,b \in \mathbb{Z}$.
 We use the factorization of $$a^6-b^6=(a^3+b^3)(a^3-b^3)=(a-b)(a+b)(a^2+ab+b^2)(a^2-ab+b^2).$$
 
Though  several solutions are obtained by combining factors,  we give four solutions as follows.

$\textbf{Case 1:} \quad W+Z=(a^3+b^3)t\ and\ W-Z=\cfrac{a^3-b^3}{t}.$
\vskip\baselineskip

We have $(W,Z)=\Bigl( \cfrac{-b^3+t^2a^3+t^2b^3+a^3}{2t},\ \cfrac{b^3+t^2a^3+t^2b^3-a^3}{2t}\Bigr)$ where $t \in \mathbb{Z}$.
\vskip\baselineskip

After canceling the denominators, we get 

$$(X,Y,Z,W)=(2at,\ 2bt,\ 4t^2(b^3+t^2a^3+t^2b^3-a^3),\ 4(-b^3+t^2a^3+t^2b^3+a^3)t^2).$$
\vskip2\baselineskip

$\textbf{Case 2:}\quad W+Z=(a^2+ab+b^2)(a^2-ab+b^2)t\ and\ W-Z=\cfrac{(a-b)(a+b)}{t}. $
\vskip\baselineskip

We have $(W,Z)=\Bigl( \cfrac{a^2+t^2a^4+t^2a^2b^2+t^2b^4-b^2}{2t},\ \cfrac{-a^2+t^2a^4+t^2a^2b^2+t^2b^4+b^2}{2t}\Bigr)$ where $t \in \mathbb{Z}$.

$$(X,Y,Z,W)=(2at,\ 2bt,\ 4(-a^2+t^2a^4+t^2a^2b^2+t^2b^4+b^2)t^2,\ 4(a^2+t^2a^4+t^2a^2b^2+t^2b^4-b^2)t^2).$$
\vskip2\baselineskip

$\textbf{Case 3:}\quad W+Z=(a^2-ab+b^2)(a-b)(a+b)t\ and\ W-Z=\cfrac{a^2+ab+b^2}{t}. $
\vskip\baselineskip

We have $(W,Z)=\Bigl( \cfrac{b^2-t^2b^4+t^2a^4-t^2a^3b+t^2ab^3+a^2+ab}{2t},\ \cfrac{-t^2b^4+t^2a^4-t^2a^3b+t^2ab^3-b^2-a^2-ab}{2t}\Bigr)$ where $t \in \mathbb{Z}$.

$$(X,Y,Z,W)=(2at,\ 2bt,\ 4(-t^2b^4+t^2a^4-t^2a^3b+t^2ab^3-b^2-a^2-ab)t^2,\ 4(b^2-t^2b^4+t^2a^4-t^2a^3b+t^2ab^3+a^2+ab)t^2).$$
\vskip2\baselineskip

$\textbf{Case 4:}\quad W+Z=(a^6-b^6)t\ and\ W-Z= \cfrac{1}{t}. $
\vskip\baselineskip

We have $(W,Z)=\Bigl( \cfrac{t^2a^6-t^2b^6+1}{2t},\ \cfrac{t^2a^6-t^2b^6-1}{2t}\Bigr)$ where $t \in \mathbb{Z}$.

$$(X,Y,Z,W)=(2at,\ 2bt,\ 4(t^2a^6-t^2b^6-1)t^2,\ 4(t^2a^6-t^2b^6+1)t^2).$$
\vskip2\baselineskip

2.2. Method-2
\vskip\baselineskip

 Let $X=a ,\ Y=b,\ W=pt+a^3,\ Z=t+b^3$ where $a,b,p,t \in \mathbb{Z}$.
\vskip\baselineskip

We obtain 
$$(-1+p^2)t+(2a^3p-2b^3)=0.$$

Hence, put $t = \cfrac{2(-a^3p+b^3)}{-1+p^2}$,
\vskip\baselineskip

then we acquire a parametric solution as follows.

$$(X,Y,Z,W)=(a(p^2-1),\ b(p^2-1),\ (p^2-1)^2(2a^3p-b^3-b^3p^2),\ (p^2-1)^2(a^3p^2-2pb^3+a^3)).$$
\vskip2\baselineskip

2.3. Method-3
\vskip2\baselineskip

Let $X=x+u,\ Y=x-u,\ W=t^2+x,\ Z=t^2-x$ where $t,u,x \in \mathbb{Z}$.

We get
\begin{equation}
12ux^4+40u^3x^2+(-4t^2+12u^5)=0.
\end{equation}

Hence, we have

$$ t^2 = 3ux^4+10u^3x^2+3u^5.$$

Take $u=19$, we get

$$Q: t^2 = 57x^4+68590x^2+7428297.$$

Quartic equation is birationally equivalent to an elliptic curve.

$$E: N^2= M^3 + M^2 - 1564M -18304.$$

We know E has rank $1$ and generator $P(M,N)=(-\cfrac{4973}{289},\ -\cfrac{303240}{4913})$ using software \textbf{SAGE}\cite{g}.
\vskip\baselineskip

The point $P(M,N)$ is of infinite order,  and the multiples $mP, m = 2, 3, \cdots$ give infinitely many points on $E$.
\vskip\baselineskip

Doubling the point $P(M,N)$ by group law using \textbf{SAGE}, we obtain $$2P(M,N)=(\frac{54247449481}{815673600},\ -\frac{9780470064206171}{23295638016000}). $$

Thus, we can obtain a quartic point corresponding to $P(M,N)$ $$2Q(x,t) = \Bigl(-\frac{14407111}{765404},\ \frac{3653190399329653}{585843283216} \Bigr).$$
\vskip\baselineskip

We give two solutions of equation $(2.1)$ using $Q$ and $2Q$.
\vskip\baselineskip

\textbf{Case Q :}

$$(X,Y,Z,W)=(57,\ 95,\ 2305248245,\ 2305248093).$$

\textbf{Case 2Q :}
\begin{align*}
&X = 51880996630,\\
&Y = 11079141384474,\\
&Z = 977318891095178497546429988862159754468626,\\
&W = 977317944922183691537359421841783861640210.
\end{align*}
\vskip3\baselineskip

\section{ Solving the diophantine equation $X^6-Y^6=W^3-Z^3$}
\vskip\baselineskip

\begin{equation}
X^6-Y^6=W^3-Z^3 \label{eq1}
\end{equation}
We give two methods to prove that there are infinitely many integer solutions of $X^6 - Y^6 = W^3 - Z^3$.
\vskip2\baselineskip

3.1. Method-1
\vskip2\baselineskip

Let $X=a,\ Y=b,\ Z=pt+a^2,\ W=t+b^2$ with $a,b,p,t \in \mathbb{Z}$.
\vskip\baselineskip

We obtain $$(p^3-1)t^2+(3a^2p^2-3b^2)t+(3a^4p-3b^4)=0.$$

Hence, take $p=\cfrac{b^4}{a^4}$ and $t= -\cfrac{3b^2a^6}{a^6+b^6}$
\vskip\baselineskip

 then we acquire a parametric solution as follows.
\vskip\baselineskip
\begin{align*}
&X = a(a^2+b^2)(a^4-a^2b^2+b^4),\\
&Y = b(a^2+b^2)(a^4-a^2b^2+b^4),\\
&Z = (a^2+b^2)(a^4-a^2b^2+b^4)b^2(b^6-2a^6),\\
&W = (a^2+b^2)(a^4-a^2b^2+b^4)a^2(a^6-2b^6).
\end{align*}
\vskip3\baselineskip

\vskip2\baselineskip
3.2. Method-2
\vskip2\baselineskip
Let $X=x+u,\ Y=x-u,\ W=t+x,\ Z=t-x$ with $t,u,x \in \mathbb{Z}$.  
\vskip\baselineskip

We have 

$$-12ux^4+(2-40u^3)x^2+(6t^2-12u^5)=0.$$

Put $y=3t$,\ hence we consider the quartic equation

$$Q: y^2 = 18ux^4+(-3+60u^3)x^2+18u^5.$$

Take $u=2$, then quartic equation is birationally equivalent to an elliptic curve.

$$E: N^2 = M^3 +M^2 -122M -444.$$

We know E has rank $1$ and  generator $P(M,N)=(17,\ 44)$ using software \textbf{SAGE}\cite{g}.
\vskip\baselineskip

The point $P(M,N)$ is of infinite order,  and the multiples $mP, m = 2, 3, \cdots$ give infinitely many points on $E$.
\vskip\baselineskip

Doubling the point $P(M,N)$ by group law using \textbf{SAGE}, we obtain $$2P(M,N)=(21,\ -93). $$

Thus, we can obtain a quartic point corresponding to $P(M,N)$ $$2Q(x,y) = \Bigl( -\frac{44}{15},\ \frac{6428}{75}\Bigr).$$
\vskip\baselineskip

We give two solutions of (3.1) using $Q$ and $2Q$.
\vskip\baselineskip

\textbf{Case Q :}

$$(X,Y,Z,W)=(6,\ 2,\ 40,\ 48).$$

\textbf{Case 2Q :}

$$(X,Y,Z,W)=(210,\ 1110,\ 1594800,\ 1297800).$$
\vskip3\baselineskip

\section{ Solving the diophantine equation $X^6 - Y^6 = R^4 - S^4$}
\begin{equation}\label{eq1}
X^6 - Y^6 = Z^4 - W^4. 
\end{equation}
We give two methods to prove that there are infinitely many integer solutions of $X^6 - Y^6 = R^4 - S^4$.
\subsection{ Method-1}

Let's consider the parametric solution of $X^6 - Y^6 = Z^2 - W^2$ where $Z=R^2$ and $W=S^2$.
Taking $X=a$ and $Y=b$. We factorize the above equation, then we get
$$(Z+W)(Z-W)=a^6-b^6=(a^3-b^3)(a^3+b^3)=(a-b)(a^2+ab+b^2)(a+b)(a^2-ab+b^2).$$
We choose $Z+W$ and $Z-W$ from the factor set $\{a-b, a^2+ab+b^2, a+b, a^2-ab+b^2\}$.
For instance, let's consider the simultaneous equations as follows.
\begin{align*}
Z + W &=(a^3+b^3)t,  \\
Z - W &=\frac{a^3-b^3}{t}. 
\end{align*}
Thus, we obtain a parametric solution of $X^6 - Y^6 = Z^2 - W^2$.
$$(X,Y,Z,W)=(2at,\ 2bt,\ 4(-b^3+t^2a^3+t^2b^3+a^3)t^2,\ 4t^2(b^3+t^2a^3+t^2b^3-a^3)).$$
We use the parametric solution of $X^6 - Y^6 = Z^2 - W^2$ as follows.
Let's consider the simultaneous equations \eqref{eq2},\eqref{eq3}.
\begin{align}
u^2 = (a^3+b^3)t^2-b^3+a^3, \label{eq2} \\
v^2 = (a^3+b^3)t^2+b^3-a^3. \label{eq3}
\end{align}
If we can find the rational solutions of equations $(4.2)$ and $(4.3)$, we can obtain the rational solutions of equation $(4.1)$. Take $(a,b)=(2,\ 1)$ then we obtain
\begin{align}
u^2 = 9t^2+7, \label{eq4} \\
v^2 = 9t^2-7. \label{eq5}
\end{align}
Let's parameterize equation \eqref{eq4}  with $(t,u)=(1,\ 4)$ and substitute $t=\cfrac{k^2-8k+9}{k^2-9}$ to \eqref{eq5}, we obtain
$$v^2 = \frac{2(k^4+432k^2-72k^3+81-648k)}{(k^2-9)^2}.$$
Hence, we consider the quartic equation
$$Q: V^2 = 2k^4-144k^3+864k^2-1296k+162.$$
The quartic equation is birationally equivalent to an elliptic curve.
$$E: y^2= x^3-49x.$$
We know E has rank $1$ and generator $P(x,y)=(\cfrac{112}{9},\  \cfrac{-980}{27})$ using \textbf{SAGE}\cite{g}.
Since the rank of E is equal to one, this implies that there are infinitely many rational points.
Using the group law, we present an example solution for equation \eqref{eq1}.
For instance, doubling the point $P$, we get a quartic point $$2Q(k,V)=\Bigl(\cfrac{508773}{142471},\ \cfrac{-362848187502}{20297985841}\Bigr).$$ 
Hence, we obtain $$(t,u,v)=\Bigl(\frac{-23058557761}{12694682160},\ \frac{25632788161}{4231560720}, \frac{20158232639}{4231560720}\Bigr).$$
Finally, we get 
\begin{align*}
X &= 5241494398998670451811238179043229184000,\\
Y &= 2620747199499335225905619089521614592000,\\
R &= 426395698533213485550190666430039500411758184842911016960000,\\
S &= 335327691756108236733334752764166296394632689254867671040000.
\end{align*}

\subsection{ Method-2}

Let $X=at,\ Y=bt,\ R=ct,\ U=\frac{S}{t}$, then \eqref{eq1} becomes 
\begin{equation}
(a^6-b^6)t^2 = c^4-U^4. \label{eq6}
\end{equation}
Put $a=2,\ b=1,\ c=4,\ V=(a^6-b^6)t$, then \eqref{eq6} becomes 
$$Q: V^2 = -63U^4 + 16128.$$
The quartic equation is birationally equivalent to an elliptic curve.
$$E: y^2= x^3\ + \ 196x.$$
We know E has rank $1$ and generator $P(x,y)=(2,\ 20)$ using \textbf{SAGE}.
Since the rank of E is equal to one, this implies that there are infinitely many rational points.
Using the group law, we present an example solution for equation \eqref{eq1}.
For instance, doubling the point $P$, we get a quartic point $$2Q(U,V)=\Bigl(\cfrac{452}{463},\ \cfrac{27175680}{214369}\Bigr).$$ 
Hence, we obtain
$$t=\cfrac{V}{a^6-b^6}=\cfrac{431360}{214369}.$$
Finally, we get 

\begin{align*}
X &= 184940423680,\\
Y &= 92470211840,\\
R &= 79290987367715840,\\
S &= 19351796053027840.
\end{align*}
\vskip3\baselineskip

\section{Concluding Remarks}
\vskip\baselineskip

We proved that there are infinitely many integer solutions of $x^n+y^n=z^n+w^n, n=2,3,4$.

However, finding the primitive solutions of $x^n+y^n=z^n+w^n, n=2,3,4$ remains an open problem.

\end{document}